\documentclass[12pt]{article}
\usepackage{amsmath}
\usepackage{amscd}
\usepackage{amssymb}
\usepackage{array}

\newtheorem{proposition}{Proposition}[section]
\newtheorem{lemma}{Lemma}[section]
\newtheorem{remark}{Remark}[section]

\begin{document}

\newcommand{\R}{\mathbb{R}}
\newcommand{\C}{\mathbb{C}}

\newcommand{\fg}{\mathfrak{g}}
\newcommand{\fh}{\mathfrak{h}}

\newcommand{\B}{\mathcal{B}}

\newcommand{\A}{\mathcal{A}}
\newcommand{\F}{\mathcal{F}}
\newcommand{\I}{\mathcal{I}}
\newcommand{\J}{\mathcal{J}}
\newcommand{\Ha}{\mathcal{H}}
\newcommand{\Oc}{\mathcal{O}}
\newcommand{\m}{\mathbf{m}}
\newcommand{\La}{\mathcal{L}}
\newcommand{\N}{\mathcal{N}_{\Theta}}
\newcommand{\U}{\mathcal{U}}
\newcommand{\de}{\delta}
\newcommand{\Z}{\mathbb{Z}}
\newcommand{\sym}{\mbox{Sym}}
\newcommand{\Gc}{G_{\C}}
\newcommand{\fgc}{\fg_{\C}}
\newcommand{\la}{\lambda}
\newcommand{\Ss}{\mathcal{S}}
\vskip 2cm
  \centerline{\LARGE \bf A Comparison Between Star Products on }
  \bigskip
 % \centerline{ \LARGE\bf  }
%\bigskip
\centerline{\LARGE\bf  Regular Orbits of Compact  Lie Groups}

\vskip 2.5cm

\centerline{R. Fioresi$^\ast$\footnote{Investigation supported by
the University of Bologna, funds for selected research topics.}
 and M. A. Lled\'o$^\dagger$}

\bigskip

\centerline{\it $^\ast$Dipartimento di Matematica, Universit\`a di
Bologna }
 \centerline{\it Piazza di Porta S. Donato, 5.}
 \centerline{\it40126 Bologna. Italy.}
\centerline{{\footnotesize e-mail: fioresi@dm.UniBo.it}}

\bigskip

\centerline{\it $^\dagger$ Dipartimento di Fisica, Politecnico di
Torino,} \centerline{\it Corso Duca degli Abruzzi 24, I-10129
Torino, Italy, and} \centerline{\it INFN, Sezione di Torino,
Italy.} \centerline{{\footnotesize e-mail:
lledo@athena.polito.it}}
 \vskip 2cm

\begin{abstract}

In this paper an algebraic  star product  and differential one
defined on  a regular coadjoint orbit of a compact semisimple
group are compared. It is proven that there is an injective algebra
homomorphism between the algebra of polynomials with the algebraic
star product and the algebra of differential functions with the
differential star product structure.
\end{abstract}

\vfill\eject

\section{Introduction}

The problem of  classification of differential star products on a
general Poisson manifold was solved in Ref.\cite{ko}. The
existence  of star products on  symplectic manifolds was already
proven in Ref.\cite{omy, dl} and, using a different technique, a
construction of a star product and a  classification of all star
products on a symplectic manifold was given in Ref.\cite{fe,de}.
For other special cases, as for regular manifolds, a proof of
existence of tangential star products was known (see Ref.
\cite{ma}).

To motivate our discussion, let us consider the Heisenberg group
$H =\R^3$ with multiplication $$(a_1,b_1,c_1)\cdot(a_2,b_2,c_2)=
(a_1+a_2,b_1+b_2,c_1+c_2+a_1b_2).$$ Its Lie algebra is
$\fh=\R^3={\rm span}\{Q,P,E'=-iE\}$ with commutation rules
$$[Q,P]=-iE\qquad \hbox{(the rest trivial)}.$$ The coadjoint
orbits of $H$ are the planes $c=\mathrm{constant}\neq 0$ (regular
orbits) and the points $(a,b,0)$. One way of obtaining the
Moyal-Weyl product on $\R^2$ is considering the Weyl map or
symmetrizer in the enveloping algebra of $\fh$,
\begin{eqnarray}\mathrm{Sym}:\mathrm{Pol}[\fh^*]&\rightarrow& U(\fh)\nonumber\\x_1x_2\cdots
x_k&\mapsto&
\frac 1{k!}\sum_{\sigma\in S_k}X_{\sigma(1)}\cdots
X_{\sigma(k)},\label{heis}
\end{eqnarray}
where $x_i$ are the coordinates on $\fh^*$ on the basis dual to
$\{X_i\}$ in $\fh$. By multiplying the commutation rules by a
formal  parameter $h$ we obtain the following star product on
$\mathrm{Pol}[\fh^*][[h]]$ \begin{equation}f\star
g=\mathrm{Sym}^{-1}(\mathrm{Sym}(f)\mathrm{Sym}(g)).\label{heis2}\end{equation}
This star product is {\it differential}, so it can be extended to
$C^\infty(\fh^*)$, it is {\it tangential}, so it can be restricted
to the orbits and it is {\it algebraic}, that is, it is closed
(and convergent) on polynomials.

If instead of the Heisenberg group we take another group, say
SU(2), we can define a star product using the symmetrizer as is
(\ref{heis}). The resulting star product is algebraic and
differential but it is not tangent to the coadjoint orbits, so it
does not define a star product on them, in this case the spheres.

Two different approaches can be taken at this point. One is to
look for a differential star product on the sphere in the spirit
of Refs. \cite{omy,dl,fe,ko}. The resulting star product is
neither algebraic nor appears related to the product on the
enveloping algebra.  The other approach insists on using the
product in the enveloping algebra. The consequence is that
differentiability is lost. This kind of star products have been
considered in Refs. \cite{cg, fl, fll, ll} and in particular, in
Refs. \cite{fl,fll} it was proven that a non differential star
product on coadjoint orbits of SU(2) corresponds to the standard
quantization of angular momentum. It seems then unavoidable to
look to a wider class of star products than the differential ones.
In particular, one cannot immediately assume that the canonical
quantization given by Kontsevich's theorem \cite{ko} is the one
relevant for physics in all cases.

 The problem of existence and classification
of algebraic star products on algebraic Poisson varieties appears
as a separate problem, mathematically interesting in itself, which
has been recently studied in Ref.\cite{ko2}. From the physical
point of view it is of interest since the algebra of a physical
quantum system may have a non differential star product, as in the
case of the angular momentum and its standard quantization.

Our purpose here is  to compare the deformations obtained by
algebraic \cite{fl,fll} and by differential methods on regular
coadjoint orbits of compact semisimple Lie groups. We want to
establish if there is some kind of equivalence among these
different star products.  We work with a family of algebraic star
products, not all isomorphic, and we relate them to the
differential star product given by Kontsevich's theorem or
Fedosovs's construction \cite{fe}. Our result is that one of the
algebraic star products can be injected homomorphically into the
differential one.

The organization of the paper is as  follows. In section 2 we
recall known facts concerning coadjoint orbits of a semisimple
compact group $G$ and its complexification $\Gc$. In section 3 we
introduce different star products on a fixed regular coadjoint
orbit $\Theta$  and on a tubular neighborhood of the orbit $\N$.
 and we prove that two different star products
on $\N$, one tangential $\star_T$ and one not tangential
$\star_{S\N}$, are equivalent. In section 4 we show our main
result, that there is an injective homomorphism between an
algebraic star product $\star_{P\Theta}$ and a differential one
$\star_{T\Theta}$ on the orbit $\Theta$. The algebraic star
product belongs to the family constructed in \cite{fl}, while the
differential one is obtained by gluing tangential star products
defined on open sets
 of  $\N$, computed with Kontsevich's formula \cite{ko}. In
 Appendix A we give for completeness some standard definitions and
 results on star products and deformations. In Appendix B we give an explicit
  formula for the gluing of star
products given in open sets and satisfying a compatibility
condition, in terms of a partition of unity.
\section{Coadjoint orbits of semisimple Lie groups}
Let $G$ be a compact semisimple group of dimension $n$ and rank
$m$ and $\fg$ its Lie algebra. Let $\fg^*$ be the dual of $\fg$.
On $C^\infty(\fg^*)$ we have the Kirillov Poisson structure: $$
\{f_1,f_2\}(\lambda)=\langle[(df_1)_{\lambda},(df_2)_{\lambda}],\lambda\rangle,
\qquad f_1,f_2 \in C^{\infty}(\fg^*), \quad \lambda \in \fg^*. $$
$(df)_{\lambda}:\fg^*\rightarrow \R$ can be considered as an
element of $\fg$, and $[\; ,\,]$ is the Lie bracket on $\fg$. Let
$\{X_1\dots X_n\}$ be  is a basis of $\fg$ and $\{x^1,\dots x^m\}$
the
 coordinates on $\fg^*$ in the dual basis. We have that
$$ \{f_1,f_2\}(x^1,\dots x^n)=\sum_{ijk}c_{ij}^k x^k\frac{\partial
f_1} {\partial x^i }\frac{\partial f_2}{\partial x^j }, $$ where
$c_{ij}^k$ are the structure constants of $\fg$, that is
$[X_i,X_j]=\sum_kc_{ij}^kX_k$.

$\fg^*$ is an algebraic Poisson manifold since   the ring of
polynomials $\R[\fg^*]$, is closed under the  Poisson bracket.

The Kirillov Poisson structure is neither symplectic nor regular.
The symplectic leaves are the orbits of the coadjoint action of
$G$ on $\fg^*$,
 $$ \langle{\rm{Ad}}^*(g)\lambda ,
Y\rangle=\langle\lambda,{\rm{Ad}}(g^{-1})Y\rangle \quad \forall\;
g\in G,\quad \lambda\in \fg^*,\quad Y\in \fg. $$ We denote by
$\Theta_\lambda$ the orbit of an element $\lambda\in\fg^*$ under
the coadjoint action.

Let $\Gc$ be the complexification of $G$ and let $\fgc$ be its Lie
algebra.  Let $\Theta_{\la \C}$ be the coadjoint orbit of
$\lambda\in \fg^*$ in $\fgc^*$ under the action of $\Gc$.
$\Theta_{\la \C}$ is an algebraic variety defined over $\R$ and $
\Theta_\lambda=\Theta_{\lambda\C} \cap \fg^*$.

Let $\C[\fg^*]$ be  the ring of polynomials on $\fgc^*$. We denote
by $ {\rm{Inv}}(\fgc^*)$ the subalgebra of  polynomials invariant
under the coadjoint action. It  is generated by homogeneous
polynomials,
  $p_i, \, i=1,\dots
m$, (Chevalley generators). We have that

$$ {\rm{Inv}}(\fg^*)= {\rm{Inv}}(\fgc^*)\cap \R[\fg^*].$$

 If $\la$
is regular,
 the ideal of $\Theta_{\la\C}$ is given by
\cite{kos}, $ \I_{0\C}=(p_i-c_i, \;i=1 \dots m), \; c_i \in \R,$
and  the polynomials on $\Theta_{\la\C}$  by $
\C[\Theta_{\la\C}]=\C[\fg^*]/\I_{0\C}$. For the real forms the
ideal of $\Theta_\lambda$ is $\I_0= \I_{0\C}\cap\R[\fg^*]$, with
the same generators than the complex one and
$\R[\Theta_{\la}]=\R[\fg^*]/\I_0=\C[\Theta_{\la\C}]\cap
\R[\fg^*]$.

\section{Star products on a regular coadjoint orbit. \label{cosp}}

In this section we will consider complex star products which are
deformations of the complexification of a real Poisson algebra. We
want to describe different star products \cite{fll} that will be
later compared.

From now on we fix a regular coadjoint orbit $\Theta$ in $\fg^*$.
We will consider $\fg_h$ the Lie algebra over $\C[[h]]$ obtained
by multiplying the structure constants of $\fg_\C$ by a formal
parameter $h$. $U_h$ is its enveloping algebra.

\paragraph{The star products $\star_S$ and $\star_{S\N}$\\}

It is well known that $U_{h}$ is a formal deformation of
$\C[\fg^*]$. In Ref. \cite{ko} it was shown that this deformation
is isomorphic to the star product canonically associated  to the
Kirillov Poisson structure. Moreover, since the linear coordinates
on $\fg^*$ are global, one can compute a star product using
Kontsevich's universal formula.

  The symmetrizer
Sym (\ref{heis}) (that can be defined in the same way for any Lie
algebra) defines through (\ref{heis2}) a differential and
algebraic star product on $\fg^*$ that we denote by $\star_S$. Any
other isomorphism that is the identity modulo $h$ could be chosen
in the place of $\sym$. All the star products constructed in this
way are equivalent to the one obtained with Kontsevich's explicit
formula. All of them are algebraic and differential, but none of
them is tangential to all the orbits \cite{cgr}.

\bigskip

Since a differential star product tangential to all orbits cannot
exist in the whole $\fg^*$ (see appendix B), we have to look for a
smaller space. We consider  a regular orbit $\Theta$ and a
regularly foliated neighborhood of the orbit, a tubular
neighborhood $\N\simeq\Theta\times \R^m$, where the global
coordinates in $\R^m$ are the invariant polynomials $p_i,\,
i=1,\dots m$. Since $\star_S$ is differential, it can be
restricted to the open set $\N$. We will denote that restriction
by $\star_{S\N}$. $\star_{S\N}$ is a differential star product
belonging to the canonical equivalence class associated to the
Kirillov Poisson structure restricted to  $\N$.

Since $\N$ is a regular Poisson manifold, we know that a
tangential star product (with respect to the symplectic leaves) exists \cite{ma}. We want to prove that
there exists a tangential star product on $\N$ equivalent to
$\star_{S\N}$.

\paragraph{The star products $\star_T$ and $\star_{T\Theta}$ \\}

We want to define a tangential star product $\star_T$ on $\N$ and
its restriction $\star_{T\Theta}$ to the regular orbit $\Theta$.
We will use the gluing of star products computed in the Appendix B
in terms of a partition of unity.

Let $\U=\{U_r,\, r\in J\}$, $J$ a set of
indices,  be a good covering of $\N$ with Darboux charts.
The coordinates in an open set $U_r$ are
\begin{eqnarray*}&&\varphi_r:U_r\longrightarrow \R^n \qquad {\rm with}\\
&&\varphi_r=(\theta_r,\pi_r,p)= (\theta_r^1, \dots
,\theta_r^{(n-m)/2}, \pi_r^1, \dots , \pi_r^{(n-m)/2},p_1,
\dots,p_m ), \\  &&\{\theta_r^\alpha,\pi_r^\beta\}=
\delta^{\alpha\beta}, \qquad\{\theta_r^\alpha,p_i\}=0,\quad
\{p_i,\pi_r^\beta\}=0.\end{eqnarray*} The invariant polynomials
$p_i$ are global coordinates, so  $U_r\simeq \hat U_r\times \R^m$
and $\{(\hat U_r, (\theta_r, \pi_r))\}_{r \in J}$ is an atlas of
$\Theta$, with $\{\hat U_r,(\theta_r, \pi_r), \, r\in J\}$ the
symplectic charts.

We can now apply Kontsevich's formula in a coordinate patch $U_r$,
using the  Darboux coordinates $\varphi_r$. We denote this star
product by $\star^K_r$. It is a tangential star product.  If
$\star_r$ denotes the restriction of $\star_{S\N}$ to $U_r$, then
$\star_r$ and $\star^K_r$ are equivalent. We will denote by
 $$
R_r:(C^\infty(U_r)[[h]],\star_r) \longrightarrow
(C^\infty(U_r)[[h]],\star_r^K) $$ the isomorphism $$ R_r(f\star_r
g)=R_r(f)\star^K_rR_r(g), \qquad  R_r={\rm Id} +\sum_{i=1}^\infty
h^nR_r^n. $$

In the intersection $U_{rs}=U_r\cap U_s$ one has that $\star^K_r$
and $\star^K_s$ are equivalent as in (\ref{intersection}) of
Appendix B with
\begin{equation}T_{rs}= R_r\circ R_s^{-1}.\label{transition}
\end{equation}
We have the following

\begin{proposition}
Let $\N$ and  $\U$ be the tubular neighborhood of the orbit
$\Theta$ and the
 covering of $\N$ defined above. Let $\F_S$ be the
sheaf of star products defined by $\star_{S\N}$ and $\star_r^K$ the
star product obtained via Kontsevich formula in $U_r \in \U$.

The assignment
$$U_r\mapsto (C^\infty(U_r), \star^K_r) \qquad \forall \;U_r\in
\U$$ is a sheaf of star products isomorphic to $\F_S$. There is a
star product $\star_T$ on $\N$
that is tangential and gauge equivalent
to $\star_{S\N}$. \label{equiv}
\end{proposition}

{\it Proof.} It is immediate that the transition functions
(\ref{transition}) satisfy the conditions (\ref{cocyco}) of
Appendix B, so we have a sheaf of star products that we will
denote by $\F_T$. The isomorphisms $R_r$ give the isomorphism  of
sheaves among $\F_S$ and $\F_T$.

Given a partition of unity subordinated to $\U$ one can use the
method of Appendix B to construct a global star product. From the
explicit formula (\ref{gsp}),  one can see that it is a tangential
star product.   \hfill$\blacksquare$

\medskip

The restriction of $\star_T$ to
the orbit will be denoted by $\star_{T\Theta}$.

\paragraph{The star products $\star_P$ and $\star_{P\Theta}$ \\}

We want to define an algebraic star product $\star_P$ on $\fg^*$
and its restriction to the orbit $\Theta$, the algebraic star
product $\star_{P\Theta}$.

We consider the ideal in $U_{h}$ $$\I_h=(P_i-c_i(h), \, i=1,\dots
m),$$ where $P_i=\sym(p_i)$ and $c_i(h)\in \C[[h]]$ with
$c_i(0)=c_i^0$. It was proven in Ref.\cite{fl} that $U_{h}/\I_{h}$
is a deformation quantization of $\C[\Theta]=\C[\fg^*]/\I_0$ where
$$\I_0=(p_i-c_i^0, \, i=1,\dots m)$$ is the ideal of a regular
orbit $\Theta$. Further properties of this deformation where
studied in Ref.\cite{fll}. The generalization of this construction
to non regular orbits was done in Ref.\cite{ll}.

A star product associated  to this deformation can be constructed
by giving a $\C[[h]]$-module isomorphism: $$
\psi:\C[\fg^*][[h]]\longrightarrow U_{h} $$ that maps the ideal
$\I_0$ isomorphically onto $\I_{h}$.  One way of choosing this map
(but not the only one) is by using the decomposition of
$\C[\fg^*]$ in terms of invariant and harmonic polynomials
\cite{kos} $$ \C[\fg^*] \cong {\rm Inv }(\fgc^*) \otimes \Ha .$$
The harmonic polynomials $\Ha$ are in one to one correspondence
with $\C[\Theta]$ and  we have a monomial basis
$\B=\{x_{i_1}\dots x_{i_k}, \; (i_1,\dots i_k)\in I \}$, where $I$
is some subset of indices such that $\B$ is a basis of
$\C[\Theta]$ (see Ref.\cite{fl} for more details). We consider the
following $\C[[h]]$-module isomorphism:
\begin{eqnarray}
\psi:\C[\fg^*][[h]] & \longrightarrow & U_{h} \nonumber\\
(p_{i_1}-c_{i_1}^0) \cdots (p_{i_k}-c_{i_k}^0)\otimes
x_{j_1}\cdots x_{j_l} & \mapsto & (P_{i_1}-c_{i_1}(h)) \cdots
\label{defpsi}\\ & & (P_{i_k}-c_{i_k}(h))\otimes (X_{j_1} \cdots
X_{j_l}),\nonumber
\end{eqnarray}
with $x_{j_1}\cdots x_{j_l}\in \B$.  $\psi$ defines an algebraic
star product on $\C[\fg^*][[h]]$, that we will denote by
$\star_P$. Since $\psi$ descends to the quotient, it also defines
an algebraic star product on $\C[\Theta][[h]]$ and  we will denote
it by $\star_{P\Theta}$. The case with $c_i(h)=c_i^0$ was
considered first in Ref.\cite{cg}, where it was shown that the
star product is not differential.

\section{Comparison between $\star_{T\Theta}$ and
$\star_{P\Theta}$}

In this section we want to compare the differential
star product $\star_{T\Theta}$ and the algebraic star product
$\star_{P\Theta}$ defined
on a fixed regular coadjoint orbit $\Theta$.
We want to show that there is an injective
algebra homomorphism $$\tilde
H:(\C[\Theta][[h]],\star_{P\Theta})\longrightarrow
(C^\infty(\Theta)[[h]]_{\C},\star_{T\Theta}).$$ We will first show that
there exists an injective algebra homomorphism
\begin{equation}H:(\C[\fg^*][[h]],\star_{P})\longrightarrow
(C^\infty(\N)[[h]]_{\C},\star_{T}),\label{ih1}\end{equation} and then
we will show that it descends appropriately to the quotients as an
injective homomorphism.

In order to compare the tangential star products $\star_P$ on
$\fg^*$ (algebraic, not differential) and $\star_T$ on $\N$ (not
algebraic, differential) we will use the non tangential star
product $\star_S$ on $\fg^*$ (algebraic and  differential).

The algebraic star products $\star_P$ and $\star_S$ on $\fg^*$ are
equivalent, since they define algebra structures that are
isomorphic to $U_{h}$. The equivalence is realized by the
$\C[[h]]$-module isomorphism:
\begin{eqnarray*}
&\eta:(\C[\fg^*][[h]],\star_P) \rightarrow
(\C[\fg^*][[h]],\star_S)&\\&\eta= \sym^{-1}\circ\psi,\qquad
\eta(f\star_P g)=\eta(f)\star_S\eta(g).&
\end{eqnarray*}
By the very definition (\ref{defpsi}) $$ f\star_Pp_i=f\cdot p_i,
$$ so, since the $p_i$'s are central the ideal $\I_0=(p_i-c_i^0)$
in $\C[\fg^*][[h]]$ with respect to the commutative
 product  is equal to the
ideal with respect to the product $\star_P$,
$\I_0^{\star_P}=(p_i-c_i^0)_{\star_P}$.

The generators of the ideal are mapped as $$
\eta(p_i-c_i^0)=(\sym^{-1}\circ\psi)(p_i-c_i^0)=\sym^{-1}(P_i-c_i(h))=p_i-c_i(h),
$$ so the ideal $\I_0^{\star_P}$ is mapped isomorphically by
$\eta$ onto the ideal with respect to the product $\star_S$,
$\I_{c(h)}^{\star_S}=(p_i-c_i(h))_{\star_S}$. We note that in the
case of $\star_S$, $\I_{c(h)}$, the ideal generated by
$p_i-c_i(h)$  with respect to the commutative product does not
coincide with $\I_{c(h)}^{\star_S}$.
Notice also that one can choose the $c_i(h)$'s arbitrarily, provided
that $c_i(0)=c_i^0$.

Since $\star_S$ is differential, it is well defined on the whole
$C^\infty(\fg^*)[[h]]_{\C}$. The commutative ideal generated by
$p_i-c_i^0$ on $C^\infty(\fg^*)[[h]]_{\C}$ will be denoted by
$\hat \I_0$. More generally, we can define $\I_{c(h)}=(p_i-c_i(h))
\subset \C[\fg^*][[h]]$, $\hat \I_{c(h)}=(p_i-c_i(h)) \subset
C^\infty(\fg^*)[[h]]_{\C}$. We have that $\I_0\subset \hat\I_0$
and $\I_{c(h)}\subset \hat \I_{c(h)}$.

\medskip

Let us consider the restriction map: $$
r:C^\infty(\fg^*)[[h]]_{\C} \longrightarrow
C^\infty(\N)[[h]]_{\C}. $$ Since the commutative product and
$\star_{S}$ are both local, the restriction $r$  is an algebra
homomorphism, between $C^\infty(\fg^*)[[h]]_{\C}$ and
$C^\infty(\N)[[h]]_{\C}$ as commutative algebras, and also between
$(C^\infty(\fg^*)[[h]]_{\C},\star_S)$ and
$(C^\infty(\N)[[h]]_{\C},\star_{S\N})$.

We consider the restriction of polynomials $r(\C[\fg^*][[h]])$.
Since a polynomial is determined by its values on any open set, we
can identify via $r$ $(\C[\fg^*][[h]],\star_S)$ with a subalgebra
of $(C^\infty(\N)[[h]]_{\C},\star_{S\N})$.

On $C^\infty(\N)[[h]]_{\C}$ there is an equivalence among $\star_{S\N}$ and
$\star_T$ (proposition \ref{equiv}). We denote it by
\begin{eqnarray*}&\rho:(C^\infty(\N)[[h]]_{\C},\star_{S\N})
\longrightarrow (C^\infty(\N)[[h]]_{\C},\star_T),&\\
&\rho(f\star_Sg)=\rho(f)\star_T\rho(g),\qquad \rho={\rm Id}
+\sum_{n=1}^\infty h^n\rho_n ,&\end{eqnarray*} where $\rho_n$ are
bidifferential operators. We have given the injective
homomorphism (\ref{ih1}) by $H=\rho\circ r\circ \eta$.

We want now to show that $H(\I_0)=\hat \J_0$, where $\hat \J_0$
is the ideal with respect to $\star_T$ in $C^\infty(\N)[[h]]_{\C}$
generated by $p_i-c_i^0$.

\medskip

We want to find out how the generators $p_i-c_i(h)$ are  mapped
under $\rho$. The scalars are mapped into scalars, since the
bidifferential operators involved in the star products  $\star_T$
and $\star_{S\N}$ are null on the constants, and so are the operators
$\rho_n$.  We need to know $\rho(p_i)$.

\begin{remark} \label{center} \end{remark}
Since $\rho$ is an isomorphism of algebras and $p_i$ belongs to
the center of \break $(C^\infty(\N)[[h]], \star_S)$, $\rho(p_i)$ must
also be in the center of $(C^\infty(\N)[[h]], \star_T)$. A
function $f$ in the center of $(C^\infty(\N)[[h]], \star_T)$ is a
function depending only on the global coordinates $f(p_1, \dots
p_m)$, since the condition $$f\star_Tg-g\star_Tf=0\qquad \forall
g\in C^\infty(\N)$$ implies for the Poisson bracket
$$\{f,g\}=0\qquad \forall g\in C^\infty(\N),$$ so $\{f,\;\}$ is a
null Hamiltonian vector field and in particular does not have
components tangent to the symplectic leaves.\hfill$\blacksquare$

\begin{remark}\label{remcenter}\end{remark}
The algebra homomorphism condition determines the form of $\rho$
on the center in terms of $\rho(p_i)=p_i+a_i(p,h)$, where
$a_i(p,h)=hz_i(p,h)$. In fact, on the center we have $$\rho=
\sum_{j_1\dots j_m}a_{1\dots m}^{j_1\dots
j_m}(p,h)\frac{\partial}{\partial p_1^{j_1}}\cdots
\frac{\partial}{\partial p_m^{j_m}}.$$ $a_{1\dots m}^{0\dots 0}=1$
and the rest of coefficients are multiples of h. In particular,
the images of $p_i$ are $\rho(p_i)=p_i+a_{1\dots i\dots m}^{0\dots
1\dots 0} (p,h)=p_i+a_i(p,h)$. Using the fact that $\star_T$ is
tangential we have $$f\star_Tp_i=f\cdot p_i\qquad \forall f\in
C^\infty(\N),$$  the homomorphism condition reads
$$\rho(p_1^{i_1}\cdots p_m^{i_m})=(p_1+a_1)^{i_1}\cdots
(p_m+a_m)^{i_m}.$$ The solution of this  equation is
\begin{equation}a_{1\dots m}^{j_1\dots j_m}=\frac{1}{j_1!\cdots
j_m!}a_1^{j_1}\cdots a_m^{j_1}.\label{eqcenter}\end{equation} In
particular, $\rho$ is trivial on the center if and only if
$a_1=\dots = a_m=0$.\hfill$\blacksquare$

\medskip

By remark $6.1$
%\cite{center}
we have that
$$
H(p_i-c_i^0)=\rho(p_i-c_i(h))=p_i-a_i(p,h)-c_i(h)=p_i-c_i^0 +h(
z_i(p,h)-\Delta_i(h)),
$$ where we have  denoted
$a_i(p,h)=hz_i(p,h)$ and  $c_i(h)=c_i^0+h\Delta_i(h)$. Since
$\Delta_i(h)$ is arbitrary, we can choose it as
\begin{equation}\Delta_i(h)=z_i(c_i^0,h). \label{deltas}\end{equation}

It is not hard to see that
$$z_i(p,h)-\Delta_i(h)=\sum_{j=1}^mb_{ij}(p_j-c_j^0)\in r(\I_0),$$
and we have $$ p_i-c_i^0 +h(z_i(p,h)-\Delta_i(h))=
\sum_{j=1}^m(\delta_{ij}+hb_{ij})(p_i-c_i^0). $$ The matrix
$(\delta_{ij}+hb_{ij})$ is invertible, so the ideal generated by
$H(p_i-c_i^0)$ in $H(\C[\fg^*][[h]])$ coincides with the ideal
generated by $(p_i-c_i^0)$ in $H(\C[\fg^*][[h]])$. (For $\star_T$,
the star ideal coincides with the commutative ideal).

\medskip

In order to state the main result we need a lemma.

\begin{lemma}
Let $\J_0$ be the ideal in $(H(\C[\fg^*][[h]]),\star_T)$ generated
by $(p_i-c_i^0)$ and let $\hat \J_0$ be the ideal in
$(C^{\infty}(\N)[[h]],\star_T)$ generated by the same generators.
Then $$ \hat \J_0 \cap H(\C[\fg^*][[h]]) = \J_0 $$
\label{intideals}
\end{lemma}

{\it Proof}. Since the product $\star_T$ is tangential to the
orbits the star ideals $\J_0$ and $\hat \J_0$ coincide with the
ideals with respect to the commutative product, so we will limit
ourselves to those.

One inclusion is obvious. For the other, let $b=\sum_{r=0}^\infty
b_rh^r \in H(\C[\fg^*][[h]])$. Assume that  \begin{equation}
H(b)=\sum_{i=1}^mf^iH(p_i-c_i^0), \label{eqa} \end{equation} where
$f^i=\sum_{r=0}^\infty f^i_rh^r\in C^\infty(\N)[[h]]_{\C}$ are not
unique. We need to prove that $f^i$ can be chosen in
$H(\C[\fg^*][[h]])$. We will show that there exist $q^i
=\sum_{r=0}^{\infty} q^i_rh^r \in \C[\fg^*][[h]]$ such that
$b_r=\sum_{i=1}^m q^i_r(p_i -c_i^0)$. This clearly will be enough.

By induction on $r$. For $r=0$, we look at the order 0 in $h$ of
the equation (\ref{eqa}) (we recall that $H={\rm Id}$ mod($h$)) $$
b_0=\sum_if^i_0(p_i-c_i^0) \in \C[\fg^*] \subset
C^{\infty}(\N)_{\C} $$ It is not hard to see that $f^i_0$ can be
chosen  in $\C[\fg^*]$, so we set $q_0^i=f_0^i$.

We go to the general case. By the induction hypothesis, we assume
that we have found $q^i_0, \dots q^i_r$, with
$$
b_0+b_1h+\cdots + b_rh^r=\sum_{i=1}^m(q^i_0+q^i_1h+\cdots
+q^i_rh^r)(p_i-c_i^0).
$$
Then,
$$
H(b)-H(b_0+b_1h+ \cdots +b_{r}h^{r})=
\sum_{i=1}^m(f^i-H(q^i_0+\cdots +q^i_rh^r)H(p_i-c_i^0),
$$
so
$$
h^{r+1}H(b_{r+1}+b_{r+2}h+ \cdots)=
h^{r+1}\sum_{i=1}^m(f^i_{r+1}-\sum_{s+t=r+1}H_s(q^i_t))H(p_i-c^0_i)
\quad \hbox{mod}(h^{r+2}).
$$
Since the ring $C^{\infty}(\N)[[h]]_{\C}$
is torsion free we have:
$$
H(b_{r+1}+b_{r+2}h+ \cdots)=
\sum_{i=1}^m(f^i_{r+1}-\sum_{s+t=r+1}H_s(q^i_t))H(p_i-c^0_i) \quad
\hbox{mod}(h).
$$
Now if we look at the order 0 in $h$
$$
b_{r+1}=\sum_{i=1}^m(f^i_{r+1}-\sum_{s+t=r+1}H_s(q^i_t))(p_i-c^0_i),
$$
as in the $r=0$ case, if we set
%we can choose $f^i_{r+1}$ in such a way that
$$
q^i_{r+1}=f^i_{r+1}-\sum_{s+t=r+1}H_s(q^i_t)
$$
it is not hard to see that it can be chosen as a
polynomial, which gives us the result. \hfill$\blacksquare$

\begin{proposition}
\label{inj} Let $\Theta$ be a regular coadjoint orbit of a compact
Lie group defined by the constrains $$ p_i-c_i^0, \qquad i=1,\dots
m. $$ There is an injective  homomorphism between the algebraic
deformation of $\C[\Theta]$ defined by $U_h/\I_h$ with $\I_h$
generated by $$P_i-c_i^0+h\Delta_i(h),$$ and the differential
deformation of $C^{\infty}(\Theta)_{\C}$
$(C^\infty(\Theta)[[h]]_{\C}, \star_{T\Theta})$, which is obtained
via Kontsevich formula (see \S 5 for more details), provided the
constants $\Delta_i(h)$ are chosen as in (\ref{deltas}).
\end{proposition}

{\it Proof}.
$H$ is an algebra  isomorphism onto its image. We   have the
commutative diagram

\begin{equation*}
\begin{CD}
\C[\fg^*][[h]]@>H>>H(\C[\fg^*][[h]])\subset
C^\infty(\N)[[h]]_{\C}\\ @VV{\pi}V @VV{\pi_H}V\\
\C[\Theta][[h]]@>\tilde H>> H(\C[\fg^*][[h]])/H(\I_0)\subset
C^\infty(\Theta)[[h]]_{\C}
\end{CD}
\end{equation*}
The last inclusion follows from \ref{intideals}.

\hfill$\blacksquare$

\begin{remark}\end{remark}
We want to note that the ideal $\I_{h}$ used in the previous
proposition  is not in general the ideal used in geometric
quantization. In fact for SU(2) it was shown in Ref.\cite{fl} that
the latter is generated by $$P-l(l+\hbar), \qquad
\hbar=\frac{h}{2\pi}.$$ ($P$ is the Casimir of $\mathfrak{su}(2)$).
But by the remark \ref{remcenter}, (\ref{eqcenter}) the ideal has,
either $c_i(h)=c_i^0$ or $\Delta^i(h)$ is an infinite series in
$h$ (an exponential).  Then we have a
contradiction.\hfill$\blacksquare$

\section*{Appendix A}

In this appendix we want to give some standard definitions on
deformations and star products that have been use throughout the
text.
\bigskip

{\bf Definition A1} {\it  Let $(\A,\{,\})$ be a Poisson algebra
over $\R$. We say that the associative algebra $\A_{[h]}$ over
$\R[[h]]$ is a formal deformation of $\A$ if
\medskip

\noindent 1. There exists an isomorphism of $\R[[h]]$-modules
$\psi: \A[[h]]\longrightarrow \A_{[h]}$;

\noindent 2.  $\psi(f_1f_2)=\psi(f_1)\psi(f_2)\; {\rm
mod(h)},\quad \forall f_1, f_2 \in A[[h]]$;

\noindent 3.
$\psi(f_1)\psi(f_2)-\psi(f_2)\psi(f_1)=h\psi(\{f_1^0,f_2^0\})\;
{\rm mod}(h^2), \quad \forall f_1,f_2\in \A[[h]]$, $f_i \equiv
f_i^0$ mod($h$), $i=1,2$. }

\bigskip

If $\A_{\C}$ is the complexification of a real Poisson algebra
$\A$ we can give the definition of {\it formal deformation} of
$\A_{\C}$ by replacing $\R$ with $\C$ in the above definition.

The associative product in $\A[[h]]$ defined by:
\begin{equation}
f \star g= \psi^{-1}(\psi(f)\cdot\psi(g)), \qquad f,g\in
\A[[h]]\label{isom}
\end{equation}
is called the {\it star product on $\A[[h]]$ induced by $\psi$.}

A star product on $\A[[h]]$ can be also defined as an associative
$\R[[h]]$-linear product given by the formula: $$ f \star
g=fg+B_1(f,g)h+B_2(f,g)h^2+\dots \in \A[[h]],\quad f,g\in \A $$
where the $B_i$'s are bilinear operators. By associativity of
$\star$ one has that $\{f,g\}=B_1(f,g)-B_1(g,f)$. So this
definition is a special case of the previous one where
$\A_h=\A[[h]]$ and $\star$ is induced by $\psi={\rm Id}$.

Two star products on $\A[[h]]$, $\star$ and $\star'$  are said to
be {\it equivalent} (or {\it gauge equivalent}) if there exists
$T=\sum_{n \geq 0}h^nT_n,$ with $T_n$ linear operators on
$\A[[h]]$, $T_0={\rm Id}$ such that $$ f\star
g=T^{-1}(Tf\star'Tg). $$

If $\A \subset C^{\infty}(M)$ and the operators $B_i$'s are
bidifferential operators we say that the star product is {\it
differential}. If in addition $\A=C^{\infty}(M)$ and $M$ is a real
Poisson manifold, we will say that $\star$ is a {\it differential
star product on $M$}.

In \cite{ko} Kontsevich classifies differential star products on a
manifold $M$ up to gauge equivalence.

\bigskip

{\bf Theorem A1} (Kontsevich, \cite{ko}) {\it  The set of gauge
equivalence classes of differential star products on a smooth
manifold $M$ can be naturally identified with the set of
equivalence classes of Poisson structures depending formally on
$h$,
 $$ \alpha=h\alpha_1+h^2\alpha_2+\cdots $$
 modulo the action
of the group of formal paths in the diffeomorphism group of $M$,
starting at the identity isomorphism. }

\bigskip

In particular, for a given Poisson structure $\alpha_1$, we have
the equivalence class of differential star products associated to
$h\alpha_1$. We will say that this is the equivalence class of
star products canonically associated to the Poisson structure
$\alpha_1$.

Also, an explicit universal formula to compute the bidifferential
operators of the star product associated to any formal Poisson
structure was given in Ref.\cite{ko} in the case of an arbitrary
Poisson structure on flat space $\R^n$. The formula depends on the
coordinates chosen, but it was also proven in Ref. \cite{ko} that
the star products constructed with different choices of
coordinates are gauge equivalent.

Let $\A_{\C}=\C[M_{\C}]$ be the coordinate ring of the complex
algebraic affine variety $M_{\C}$ defined over $\R$ whose real
points are  a real algebraic Poisson variety $M$. If the $B_i$'s
are bilinear algebraic operators we will say that $\star$ is an
{\it algebraic star product on $M$}.

An example of great interest for us, of such $M$ is given by the
dual $\fg^*$ of the Lie algebra of a compact semisimple Lie group
(see section 2).

The classification of algebraic star products is still an open
problem \cite{ko2}.
\bigskip

{\bf Definition A2} {\it Let $N$ be a submanifold of the Poisson
manifold $M$ and let $\star_M$ be a star product on $M$. We say
that $\star_M$ is {\it tangential} to $N$ if for $f,g \in
C^{\infty}(N)$: $$ f \star_N g =_{def}(F \star_M G)|_N, \quad
{\rm\it  with} \quad f=F|_N, \quad g=G|_N $$ is a well defined
star product on $N$, that is, if
\begin{equation}
(F-F')|_N=(G-G')|_N=0, \quad {\rm\it then} \quad F\star_M
G|_N=F'\star_M G'|_N, \label{tangential}
\end{equation}
for $F,F',G,G'\in C^\infty(M)$. }

\bigskip

 The same definition works for algebraic Poisson varieties,
replacing the algebra of $C^{\infty}$ functions with the algebra
of polynomials on the varieties.

Given a Poisson manifold $M$, one can ask if there exists a
differential star product on $M$, that is tangential to all the
leaves of the symplectic foliation. For regular manifolds a
positive answer was found in Ref.\cite{ma}. For $M=\fg^*$ foliated
in coadjoint orbits, it was found in Ref.\cite{cgr} that there is
an obstruction to the existence. In particular,  for a semisimple
Lie algebra $\fg^*$ it is not possible to find a differential star
product on $\fg^*$ which is tangential to all coadjoint orbits.

\section*{Appendix B}

In this Appendix we want to give an explicit formula on how to
construct a global star product starting from  star products
defined on open sets of a manifold and satisfying certain
conditions (see below). We will refer to this procedure as {\it
gluing of star products}, and it  will be used  in section
\ref{cosp}.

 Let $\star$ be a differential star
product on a manifold $M$. Since the operators $B_i$ that define
$\star$ are local, there are well defined star products $\star_U$
on every open set  $U$ of $M$. We have a sheaf of algebras $\Ss$:
\begin{equation}
\Ss(U)=(\C^\infty(U)[[h]],\star_U).\label{lsheaf}
\end{equation}
which we will call {\it sheaf of star products}.

Let $M$ be a Poisson manifold and fix an open cover $\U=\{U_r\}_{r
\in J}$ where $J$ is some set of indices. Assume that in each
$U_r$ there is a differential star product $$
\star_r:C^\infty(U_r)[[h]]\otimes
C^\infty(U_r)[[h]]\longrightarrow C^\infty(U_r)[[h]]. $$ This
defines a collection of sheaves of star products
\begin{equation}
\F_r(V_r)=(C^\infty(V_r)[[h]],\star_r),\qquad V_r\subset
U_r.\label{stsheaf}
\end{equation}

It is a general fact in theory of sheaves  that if there are
isomorphisms of sheaves in the intersections \begin{eqnarray}
&&T_{sr}:\F_r(U_{rs})\longrightarrow\F_s(U_{sr}),\qquad
U_{r_1\dots r_k}=U_{r_1}\cap \cdots \cap U_{r_k} \nonumber\\
&&T_{sr}(f)\star_s
T_{sr}(g)=T_{sr}(f\star_rg)\label{intersection}\end{eqnarray} such
that the following conditions are satisfied
\begin{eqnarray}
&1.& T_{rs}=T_{sr}^{-1} \qquad\qquad\; {\rm{on}} \;
U_{sr},\nonumber
\\ &2.& T_{ts}\circ T_{sr}=T_{tr}, \qquad {\rm{on}}\quad
U_{rst},\label{cocyco}\end{eqnarray} then there exists a global
sheaf $\F$ on $M$ isomorphic to the local sheaves $\F_r$ on each
$U_r$.

If the sheaves of star products (\ref{stsheaf}) satisfy the
conditions (\ref{cocyco}) with $$T_{sr}={\rm Id}\; \quad {\rm
mod}(h),$$ then we have a global sheaf of star products on $M$.
The algebra of the global sections is $C^\infty(M)[[h]]$ together
with a star product that we will call the {\it gluing} of local
star products. We want to write an explicit formula for the star
product of global sections.

We denote  $U^{r_1\dots r_k}=U_{r_1}\cup\dots \cup U_{r_k}$. Let
us first consider the gluing on two open sets, say $U_1$ and
$U_2$, with non trivial intersection. Let $\phi_1:U_1\rightarrow
\R$, $\phi_2:U_2\rightarrow \R$ be a partition of unity of
$U^{12}$, $$\phi_1(x) +\phi_2(x)=1 \quad \forall x\in
U^{12};\qquad {\rm{supp}}(\phi_r)\subset U_r.$$ Let $f_r\in
C^\infty(U_r)[[h]]$ such that $f_s=T_{sr}f_r$ in $U_{rs}$.
 One can
define an element $f\in C^\infty(U^{rs})[[h]]$ by
$f=\phi_1f_1+\phi_2f_2$. On the intersection $U_{12}$ one has
$$f=\phi_1f_1+\phi_2T_{21}f_1=(\phi_1{\rm{Id}}+\phi_2T_{21})f_1=A_{21}f_1.$$
Notice that the operator $A_{21}={\rm{Id}}+\Oc(h)$ is invertible.
On $U^{12}$ we can define the star product
\begin{equation}
f\star g=\left\{\begin{array}{ll} (f_1\star_1g_1)(x)&{\rm{if}}\;
x\in U_1-U_{12}\\ &\\
A_{21}(A^{-1}_{21}(f)\star_1A^{-1}_{21}(g))(x)&{\rm{if}}\; x\in
U_{12}\\ & \\ (f_2\star_2g_2)(x)&{\rm{if}}\; x\in
U_2-U_{12}\end{array}\right.\label{gsp}\end{equation} It is easy
to check that the star product is smooth.

One can do the gluing interchanging $U_1$ and  $U_2$. One has that
on $U_{12}$
$$f=\phi_1f_1+\phi_2f_2=(\phi_1T_{12}+\phi_2{\rm{Id}})f_j=A_{12}f_2.$$
$A_{12}$ is also invertible and $$A_{21}=A_{12}T_{21}$$ provided
$T_{12}=T_{21}^{-1}$. One can construct a star product on $U^{12}$
using  the same procedure than in (\ref{gsp}). It is easy to check
that both star products are identical.

The procedure in (\ref{gsp}) can be generalized to an arbitrary
number of open sets. Let $\phi_i:U_i\rightarrow \R$ a partition of
unity of $M$ subordinate to the covering $\U$.  We define $f\in
C^\infty(M)$ $$f=\sum_{r\in J}\phi_rf_r, \qquad {\rm{where}}\;
f_r=T_{rs}f_s.$$ On $U_r$ $f$ becomes
$$f=(\phi_r{\rm{Id}}+\sum_s\phi_sT_{sr})f_r=A_rf_r.$$ The star
product on $U_r$ is defined as \begin{equation}f\star
g=A_r(A_r^{-1}(f)\star_rA_r^{-1}(g)).\label{mgsp}\end{equation}
Using conditions (\ref{cocyco}) one has $$A_rT_{rt}=A_t.$$ Then,
the star products
 (\ref{mgsp}) on each $U_r$ coincide in the intersections, so they
define a unique star product on $M$. The restriction  of this star
product to $U_r$ is equivalent to $\star_r$.  Also, using
different partitions of unity one obtains equivalent star
products.

\section*{Acknoledgements}

\medskip

We wish to thank A. Levrero for helpful discussions.


\begin{thebibliography}{99}
\bibitem{ko} M. Kontsevich, {\it Deformation Quantization of Poisson
Manifolds}. math.QA/9709040.

\bibitem{omy} H. Omori, Y. Maeda and A. Yoshioka, {\it Weyl
Manifolds and Deformation Quantization}. Adv. in Math. {\bf 85}
224-255 (1991).

\bibitem{dl} M. De Wilde, P. B. A. Lecomte, {\it Existence of
Star Products and of Formal Deformations  of the Poisson Lie
Algebra of Arbitrary Symplectic Manifolds}. Lett. Math. Phys.
{\bf 7} 487-49 (1983).



\bibitem{fe}  B. Fedosov,  {\it A Simple Geometric Construction of
Deformation Quantization}. J. Diff. Geom. {\bf 40} Vol. 2 213-238
(1994).

\bibitem{de} P. Deligne, {\it D\'eformations de l'Algebre des Fonctions
 d'une Vari\'et\'e Symplectique: Comparaison entre Fedosov et De Wilde, Lecomte}
  Selecta Math. New series {\bf 1} No.4 667-697 (1995).


\bibitem{ma} M. Masmoudi, {\it Tangential deformation of a Poisson bracket
 and tangential star-products on a regular
Poisson manifold.}  J.  Geom.  Phys. {\bf 9} 155-171 (1992).

\bibitem {cg} M. Cahen, S. Gutt, {\it Produits $*$ sur les
    Orbites des Groupes Semi-Simples
    de Rang 1}, C.R. Acad. Sc. Paris  {\bf 296} (1983),
    s\'erie I, 821-823; {\it An Algebraic Construction
    of $*$ Product on the Regular
        Orbits of Semisimple Lie Groups}. In  ``Gravitation
and Cosmology". Monographs and Textbooks in Physical Sciences. A
volume in honor of Ivor Robinson, Bibliopolis. Eds W. Rundler and
A. Trautman, (1987);
 {\it Non Localit\'e d'une D\'eformation Symplectique sur la
  Sph\`ere $S^2$}.  Bull. Soc. Math. Belg. {\bf 36 B} 207-221 (1987).




  \bibitem{fl} R. Fioresi and M. A. Lled\'o, {\it On the deformation
Quantization of Coadjoint Orbits of Semisimple Lie Groups}. To
appear in the Pacific J. of Math. math.QA/9906104.


\bibitem{fll}  R. Fioresi, A. Levrero and M. A. Lled\'o, {\it Algebraic and
Differential  Star Products  on  Regular Orbits of Compact Lie
Groups}. To appear in the  Pacific J. of Math. math.QA/0011172.


\bibitem{ll} M. A . Lled\'o, {\it  Deformation Quantization
of Non Regular Orbits of Compact Lie Groups}.  math.QA/0105191.

\bibitem{ko2} M. Kontsevich, {\it Deformation Quantization of algebraic
varieties}. math.AG/0106006.










\bibitem{kos} B. Kostant, {\it Lie Group Representations on Polynomial
Rings}. Am. J. Math. {\bf 85} 327 (1978).

\bibitem {cgr} M. Cahen, S. Gutt, J. Rawnsley, {\it On Tangential Star
Products for the Coadjoint Poisson Structure}. Comm. Math. Phys,
{\bf 180} 99-108 (1996).





\end{thebibliography}
\end{document}